\documentclass{amsart}
\usepackage{amssymb}
\usepackage{graphicx}
\usepackage[all]{xy}
\usepackage{hyperref}

\newtheorem{lemma}{Lemma}[section]
\newtheorem{theorem}[lemma]{Theorem}
\newtheorem{thm}{Theorem}
\newtheorem{prop}[thm]{Proposition}

\newtheorem{proposition}[lemma]{Proposition}
\newtheorem{corollary}[lemma]{Corollary}

\theoremstyle{remark}
\newtheorem{problem}[lemma]{Problem}
\newtheorem{remark}[lemma]{Remark}
\newtheorem*{examples}{Examples}

\theoremstyle{definition}

\renewcommand{\phi}{\varphi}
\newcommand{\bbA}{\mathbb{A}}
\newcommand{\bbC}{\mathbb{C}}
\newcommand{\bbN}{\mathbb{N}}

\newcommand{\bbQ}{{\mathbb{Q}}}
\newcommand{\bbR}{{\mathbb{R}}}

\newcommand{\gal}{{\rm Gal}}

\newcommand{\calN}{{\mathcal N}}
\newcommand{\Mgal}{\tilde M}
\newcommand{\Fgal}{\tilde F}

\newcommand{\Khat}{\hat K}

\newcommand{\nek}{{,\ldots,}}
\newcommand{\Char}{\textnormal{char}}

\newcommand{\bfsigma}{\mbox{\boldmath$\sigma$}}
\newcommand{\isom}{\cong}
\newcommand{\frakp}{\mathfrak{p}}
\newcommand{\frakP}{\mathfrak{P}}

\newcommand{\lcm}{\textrm{lcm}}

\begin{document}
\title{Dirichlet's Theorem for Polynomial Rings}

\author{Lior Bary-Soroker}

\address{School of Mathematical Sciences
Tel Aviv University Ramat Aviv, Tel Aviv 69978 ISRAEL}

\email{barylior@post.tau.ac.il}


\begin{abstract}
We prove the following form of Dirichlet's theorem for polynomial
rings in one indeterminate over a pseudo algebraically closed field
$F$. For all relatively prime polynomials $a(X), b(X)\in F[X]$ and
for every sufficiently large integer $n$ there exist infinitely many
polynomials $c(X)\in F[X]$ such that $a(X) + b(X)c(X)$ is
irreducible of degree $n$, provided that $F$ has a separable
extension of degree $n$.
\end{abstract}

\maketitle

\section*{Introduction}
\noindent Dirichlet's classical theorem on primes in arithmetic
progressions states that if $a,b$ are relatively prime positive
integers, then there are infinitely many $c\in\bbN$ such that $a+bc$
is a prime number. Following a suggestion of Landau, Kornblum proved
an analog of Dirichlet's theorem for the ring of polynomials $F[X]$
over a finite field $F$ \cite{Kornblum1919}. Later, Artin refined
Kornblum's result and proved that if $a(X),b(X)\in F[X]$ are
relatively prime, then for every sufficiently large integer $n$
there exists $c(X)\in F[X]$ such that $a(X)+b(X)c(X)$ is irreducible
of degree $n$ \cite[Theorem 4.8]{Rosen2002}.

To avoid repetition, we shall say that Dirichlet's theorem holds for
a polynomial ring $F[X]$ and a set of positive integers $\calN$, if
for any relatively prime polynomials $a,b\in F[X]$ there exist
$n_0>0$ and infinitely many $c\in F[X]$ such that $a+bc$ is
irreducible of degree $n$, for any $n\geq n_0$ in $\calN$.

Jarden raised the question of whether the Artin-Kornblum result can
be generalized to other fields. Of course, if $F$ is algebraically
closed, then the polynomial $a(X)+b(X)c(X)$ is reducible unless it
is of degree $1$. On the other hand, if $F$ is Hilbertian, then
there are infinitely many $\alpha\in F$ such that $a(X)+b(X)\alpha$
is irreducible in $F[X]$. To get irreducible polynomials of higher
degree in this case, one may first choose $c(X)\in F[X]$ relatively
prime to $a(X)$ and of high degree, and then find $\alpha\in F$ for
which $a(X)+b(X)c(X)\alpha$ is irreducible over $F$.

Artin's proof of the result quoted in the first paragraph is based
on a weak form of Weil's theorem on the Riemann hypothesis for
absolutely irreducible curves over finite fields. The theorem
roughly states that if a finite field $F$ is large compared to the
coefficients of the equations defining the curve, then the curve has
$F$-rational points. This makes it plausible that for infinite
fields $F$ with the latter property the ring $F[X]$ should satisfy
Dirichlet's theorem. Such fields are called PAC (Pseudo
Algebraically Closed). Explicitly, a field $F$ is {\bf PAC} if every
nonempty absolutely irreducible variety defined over $F$ has an
$F$-rational point. (See \cite[Chapter 11]{FriedJarden2005} for a
comprehensive discussion of PAC fields and \cite[Theorem
18.6.1]{FriedJarden2005} for an abundance of algebraic extensions of
countable Hilbertian fields which are PAC).

Of course, if $F$ is algebraically closed, then it is PAC, but, as
pointed out above,
Dirichlet's theorem does not hold for $F[X]$ (for any infinite $\calN$). %
Let $\calN(F)$ be the set of all positive integers $n$ such that $F$
has a separable extension of degree $n$. Our main result asserts
that Dirichlet's theorem holds for $\calN(F)$:

\begin{thm}\label{thm A}
Let $F$ be a PAC field. Then Dirichlet's theorem holds for $F[X]$
and $\calN(F)$.
\end{thm}

The proof of Theorem A uses a weak form of Hilbert's Irreducibility
Theorem that PAC fields satisfy and then argues as in the third
paragraph. Roquette was the first to observe that a PAC field which
has a rich Galois structure (namely, is $\omega$-free) is Hilbertian
\cite[Corollary 27.3.3]{FriedJarden2005}.

We elaborate Roquette's approach and show in
Corollary~\ref{cor:places} that if $F$ is PAC, $f\in F[X,Y]$ is
irreducible, and the splitting field of $f(X,Y)$ over $F(Y)$ is
regular over $F$, then, under some necessary assumptions, there are
infinitely many specializations $Y\mapsto \alpha\in F$ for which
$f(X,\alpha)$ remains irreducible over $F$.

Note that in order to get an irreducible specialization for a
polynomial $f(X,Y)$ Roquette finds a specialization that preserve
the Galois group $\gal(f(X,Y),F(Y))$. Therefore it is somewhat
unexpected that we can find an irreducible specialization even if
$\gal(f(X,Y),F(Y))$ does not occur as Galois group over $F$.

As a preparation to the use of Corollary~\ref{cor:places}, we prove
a result over an arbitrary infinite field which is interesting for
its own sake:

\begin{prop}\label{thmB}
Let $F$ be an infinite field with an algebraic closure $\Fgal$ and
let $a(X),b(X)\in F[X]$ be relatively prime polynomials. Then for
every sufficiently large positive integer $n$ there exists $c(X)\in
F[X]$ for which $f(X,Y)=a(X)+b(X)c(X)Y$ is irreducible over $F(Y)$
of degree $n$ in $X$ and $\gal(f(X,Y),\Fgal(Y))\isom S_n$.
\end{prop}

Finally, note that each infinite algebraic extension $F$ of a finite
field $K$ is PAC \cite[Corollary 11.2.4]{FriedJarden2005}. By
Theorem~\ref{thm A},  Dirichlet's theorem holds for $F[X]$ and
$\calN(F)$. This result already follows from a
quantitative form of the result of Artin-Kornblum. %
Nevertheless, our proof has the advantage that the constructions are
essentially explicit: The polynomial $c(X)$ in Theorem~\ref{thm A}
equals to the polynomial $c(X)$ appearing in Proposition~\ref{thmB}
times some factor, say $\alpha$, coming from the PACness property.
The construction in Proposition~\ref{thmB} is explicit as it uses
nothing but the Euclidean algorithm.

\textit{Notation.}  Throughout this paper we denote by $F$ an
infinite field, by $F[X]$ and $F[X,Y]$ the polynomial rings over $F$
in one and two variables, respectively, and by $\Fgal$ a fixed
algebraic closure of $F$. As mentioned earlier $\calN(F)$ denotes
the set of all positive integers such that $F$ has a separable
extension of degree $n$. If we have a Galois extension, say $K/F$,
then $\gal(K/F)$ denotes its Galois group. The absolute Galois group
of $F$ is denoted by $\gal(F)$, i.e., $\gal(F) = {\rm
Aut}(\Fgal/F)$. For a polynomial $a \in F[X]$ we write $a'$ for its
derivative. Finally, we abbreviate and say ``for large $n$'' instead
of ``$(\exists n_0\in \bbN)(\forall n>n_0)$.''

\textsc{Acknowledgment:}  I would like to thank Moshe Jarden for
raising the question that initiated this work and for many helpful
suggestions and to Peter M\"uller for his help in the second
chapter. I would also like to thank  Dan Haran, Joseph Bernstein,
and Sasha Sodin. I gratefully thank the anonymous referee for his
valuable comments which improved this work. Special thanks are
directed to my friends Dubi Kelmer and Ilya Surding and to my wife
Hamutal for carefully reading the manuscript.

This work was partially carried out while the author was at the
Max-Planck-Institut f\"ur Mathematik in Bonn. It is part of the
author's Ph.D. thesis done at Tel Aviv University, supervised by
Prof. Dan Haran.

\section{Field crossing argument}\label{sec:1}
Let $K$ be a finitely generated regular extension of a field $M$.
Suppose we have finite Galois extensions $E/K$ and $N/M$, $E$
regular over $K$, together with an embedding $\gamma \colon A \to
G$, where $A=\gal(N/M)$ and $G=\gal(E/K)$. Identify $\gal(EN/K)$
with $G\times A$  and let $\Delta = \{(\gamma(\sigma),\sigma)\in
\gal(EN/K)\mid \sigma\in A\}$.

The famous field crossing argument uses $\gamma$ to ``twist'' $E/K$
to $D/K$, where $D$ is the fixed field of $\Delta$ in $EN$. Let
$\psi$ be an unramified $M$-place of $D$ of degree $1$ over $K$.
Extend it to an $N$-place of $DN=NE$. Its restriction to $E$ is an
unramified $M$-place $\phi$ of $E$ over $K$, with decomposition
group $\gamma(A)$ and residue field $N$ \cite[proof of Lemma
24.1.1]{FriedJarden2005}. Moreover, the canonical homomorphism
$\phi^*\colon A\to G$ defined by $\phi$ is exactly $\gamma$
\cite[Remark on page 9]{FriedHaranJarden1984}.

Next assume that $N$ and $E$ are the Galois closures of some
separable extensions $N'/M$ and $E'/K$, respectively, of degree $n$.
Let $A' = \gal(N/N')$ and $G' = \gal(E/E')$ be their respective
Galois groups and let $\Sigma = A/A'$ and $\Theta = G/G'$ be the
corresponding sets of left cosets. Then $|\Sigma|=|\Theta|=n$ and
$A$ and $G$ act naturally (by left multiplication) on $\Sigma$ and
$\Theta$, respectively. A key observation is that if the above
embedding $\gamma$ respects this extra structure (i.e., there exists
an injection, and hence bijection, $\gamma^*\colon \Sigma\to \Theta$
such that $\gamma(a)\gamma^*(\sigma)=\gamma^*(a\sigma)$, for all
$a\in A$, $\sigma\in \Sigma$, or, in other words, $\gamma\colon
(A,\Sigma)\to (G,\Theta)$ is an embedding of permutation groups),
then a stronger conclusion holds:

\begin{lemma} \label{lem:primes}
In the notation and under the assumptions above, any $M$-place
$\phi$ of $E$ for which $\phi^*=\gamma$ restricts to an $M$-place
$\phi'$ of degree $n$ of $E'$, unramified over $K$ with residue
field $N'$.
\end{lemma}

\begin{proof}
Viewing $A'$ as an element in $\Sigma$ we have $\gamma^*(A') = g G'$
for some $g\in G$. Without loss of generality we may assume that
$g=1$, that is $\gamma^*(A') = G'$. (Otherwise, we replace $(\gamma,
\gamma^*)$ with $(\gamma',\gamma'^*)$, where $\gamma'(a) = g^{-1}
\gamma(a) g$ and $\gamma'^*(\sigma) = g^{-1} \gamma^*(\sigma)$,
$a\in A$ and $\sigma\in A/A'$.) As $A'$ is the stabilizer of itself
(in $A$) and $G'$ is the stabilizer of $G'$ (in $G$), we have
$\gamma(A') = G'\cap \gamma (A)$.
Let $E''$ be the decomposition field of $\psi$, i.e., the fixed
field of $\gamma(A)$. Then $\gal(E/E'E'') = G'\cap \gamma(A) =
\gamma(A')$. Therefore the residue field of $E'E''$ under $\phi$ is
$N'$, and hence, so is the residue field of $E'$. In particular, the
degree of $\phi' = \phi|_{E'}$ is $[N':M] = n$.
\end{proof}

In what follows we shall apply the previous lemma to get a weak form
of Hilbert's Irreducibility Theorem for PAC fields and, more
generally, for fields which have a PAC extension.

A polynomial $f(X,Y)\in F[X,Y]$ is called \textbf{absolutely
irreducible} if $f(X,Y)$ is irreducible over $\Fgal$. If in addition
the Galois groups of $f(X,Y)$ over $F(Y)$ and over $\Fgal(Y)$ are
equal, then $f(X,Y)$ is said to be \textbf{$X$-stable} over $F$.
Thus, a polynomial $f(X,Y)$ is $X$-stable if and only if it is
irreducible and its splitting field over $F(Y)$ is regular over $F$
\cite[Remark 16.2.2]{FriedJarden2005}.

\begin{examples}
\renewcommand{\labelenumi}{\roman{enumi}.}
\begin{enumerate}
\item
Every absolutely irreducible polynomial $f(X,Y)$ which is Galois
over $F(Y)$ is $X$-stable over $F$.
\item
Every polynomial $f(X,Y)$ of degree $n$ in $X$ with symmetric Galois
group $S_n$ over $\Fgal(Y)$ is $X$-stable.
\item (Jarden)
If an absolutely irreducible polynomial $f(X,Y)$ has a simple Galois
group $G$ over $F(Y)$, then $f$ is $X$-stable over $F$ (since the
Galois group of $f(X,Y)$ over $\Fgal(Y)$ is a nontrivial normal
subgroup of $G$).
\item
\cite{MalleMatzat1999} contains many explicit stable polynomials
over $\bbQ$ (and hence over any field of characteristic zero), e.g.,
$f(X,Y) = X^n - Y(nX-n+1)$ is $X$-stable with symmetric Galois group
over $\bbQ(Y)$ \cite[Theorem 9.4]{MalleMatzat1999}.
\end{enumerate}
\end{examples}

A field extension $M/F$ is said to be a \textbf{PAC extension} if
for every nonempty variety $V$ of dimension $r$ defined over $M$ and
for every dominating separable rational map $\phi\colon V \to
\bbA^r$ over $M$ there exists an $M$-point $a\in V(M)$ such that
$\phi(a)\in F^r$ \cite{JardenRazon1994}. PAC extensions generalize
PAC fields since $F$ is a PAC field if and only if $F/F$ is a PAC
extension. (Note that an extension $M/F$ such that $M$ is a PAC
field, need not be a PAC extension, see Section~\ref{sec:PAC} for
details.) The property of a PAC extension can be reformulated in
terms of places: $M/F$ is a PAC extension if and only if for every
separable finite extension $D/M(Y_1,\ldots,Y_r)$ with $D$ regular
over $M$, there exist infinitely many $M$-places $\phi$ of $D/M$ of
degree $1$ such that $\phi(Y_i)\in F$, $i=1,\ldots,r$.

The following proposition establishes a weak form of Hilbert's
Irreducibility Theorem for stable polynomials over a field which has
a PAC extension.

\begin{proposition}
\label{thm:places}%
Let $M/F$ be a PAC extension, let $f(X,Y)\in M[X,Y]$ be an
$X$-stable polynomial over $M$ of degree $n$ in $X$, and let $N'/M$
be a separable extension of degree $n$ with Galois closure $N$.
Consider $\gal(N/M)$ and $\gal(f(X,Y),M(Y))$ as permutation groups
of degree $n$ via the action on the cosets of $\gal(N/N')$ in
$\gal(N/M)$ and the action on the roots of $f(X,Y)$ over $M(Y)$,
respectively. If there exists an embedding of permutation groups
$\gamma\colon \gal(N/M) \to \gal(f(X,Y),M(Y))$, then there exist
infinitely many $\alpha\in F$ for which $f(X,\alpha)$ is irreducible
over $M$. Moreover, $N'$ is generated by a root of $f(X,\alpha)$
over $M$.
\end{proposition}

\begin{remark}
The assumption on $\gamma$ is necessary. Indeed, $\gal(N/M)$ is
isomorphic as a permutation group to $\gal(f(X,\alpha),M)$ which is
a subgroup of $\gal(f(X,Y),M(Y))$ (for all but a finite number of
$\alpha$'s).
\end{remark}

\begin{proof}
Let $K = M(Y)$, let $E$ be the splitting field of $f(X,Y)$ over
$M(Y)$, and let $E' = K(x)$, where $x\in E$ is a root of $f(X,Y)$.
Then $E$ is the Galois closure of $E'/K$ and $n=[E':K]$. The field
crossing argument  gives a field extension $D/K$, regular over $M$,
with the property that any $M$-place $\psi$ of degree $1$ of $D$,
unramified over $K$, yields an $M$-place $\phi$ of $E$ for which
$\phi^* = \gamma$. Now by Lemma~\ref{lem:primes} such $\phi$
restricts to an $M$-place $\phi'$ of degree $n$ of $E'$, unramified
over $K$ over the same place of $K$ with residue field $N'$. As
$M/F$ is PAC, there exist infinitely many $M$-places $\phi$ of $D$
of degree $1$ unramified over $K$ such that $\alpha = \phi(Y)\in F$.
Then the corresponding place $\phi'$ of $E'$ has residue field $N'$.
But the residue field of $E'$ is generated by a root of
$f(X,\alpha)$, so $f(X,\alpha)$ is irreducible.
\end{proof}

The last result of this section deals with the special case where
the Galois group of the stable polynomial is the symmetric group. In
this case the condition on $\gamma$ is redundant.

\begin{corollary}
\label{cor:places}%
Let $M/F$ be a PAC extension, let $f(X,Y)\in M[X,Y]$ be a polynomial
of degree $n$ in $X$, and let $N/M$ be a separable extension of
degree $n$. Assume that the Galois group of $f(X,Y)$ over $\Mgal(Y)$
is $S_n$. Then there exist infinitely many $\alpha\in F$ for which
$f(X,\alpha)$ is irreducible over $M$ and $N$ is generated by a root
of $f(X,\alpha)$ over $M$.
\end{corollary}

\section{Polynomials over infinite fields}\label{sec:2}
\subsection{Technical background and basic tools}
The following result is a special case of Gauss' Lemma.
\begin{lemma}
\label{lem1}%
A polynomial $f(X,Y) =a(X)+b(X)Y\in F[X,Y]$ is irreducible if and
only if $a(X)$ and $b(X)$ are relatively prime.
\end{lemma}

\begin{lemma}
\label{lem_separable}%
Let $a,b,c\in F[X]$ such that $\gcd(a,b) = 1$ and $c\neq 0$. Then
there exists a finite subset $S\subseteq F$ such that for each
$\alpha\in F\smallsetminus S$ the polynomials $a + \alpha b$ and $c$
are relatively prime. Moreover, if $b'\neq 0$, we may choose $S$
such that $a +\alpha b$ is also separable.
\end{lemma}

\begin{proof}
Let $S = \{-\frac{a(\gamma)}{b(\gamma)}\mid \gamma\in \Fgal,
b(\gamma)\neq 0 \mbox{ and } c(\gamma) = 0\}\cap F$. Then $a +
\alpha b$ and $c$ have no common zero in $\Fgal$, for any
$\alpha\not\in S$. Hence these polynomials are relatively prime.
Next let $d(Y)\in F[Y]$ be the discriminant of $a(X) + Y b(X)$ over
$F(Y)$, then $b'(X)\neq 0$ implies that $d(Y)\neq 0$. In this case
add all the roots of $d(Y)$ to $S$.
\end{proof}

\begin{lemma}
\label{lem_final_poly}%
Let $a,b,p_1,p_2\in F[X]$ be pairwise relatively prime polynomials
and let $\alpha_1, \alpha_2 \in F$ be distinct nonzero elements.
Then for any $n > \deg p_1 + \deg p_2$ there exists $c\in F[X]$ of
degree $n$ and separable $h_1,h_2\in F[X]$ such that $a = p_i h_i +
b c \alpha_i$ and $\gcd(h_i, a p_i) = 1$ for $i=1,2$.
\end{lemma}

\begin{proof}
Write $b_i = b \alpha_i$. Since $\gcd(p_i,b_ip_{3-i}) = 1$ for
$i=1,2$, we have
\begin{equation} \label{eq a 0}
a = p_i h_{i,0} + b_i  p_{3-i} c_i\qquad i=1,2,
\end{equation}
with $\deg c_i <\deg p_i$. For $\bar c = p_1c_2  + p_2 c_1$ and
$h_{i,1} = h_{i,0} - b_i c_{3-i}$, we have
\begin{equation}\label{eq a 1}
a = p_i h_{i,1} + b_i \bar c \qquad i=1,2.
\end{equation}
Here $h_{i,1}$ is relatively prime to $b_i$, since $a$ and $p_i$
are. Taking \eqref{eq a 0} with $i=2$ and \eqref{eq a 1} with $i=1$
modulo $p_2$, we get
\[
p_1 h_{1,1} \equiv a - b_1 \bar c \equiv a - b_1p_1 c_{2}   \equiv
b_2p_1 c_2  - b_1 c_2 p_1 \equiv b p_1  c_2 (\alpha_2 - \alpha_1)
\pmod{p_2}.
\]
Therefore $h_{1,1}$ is relatively prime to $p_2$, since $bp_1c_2$ is
(by \eqref{eq a 0} with $i=2$). Similarly, $h_{2,1}$ is relatively
prime to $p_1$.

Take $c = \bar c + p_1 p_2 s$ for some $s\in F[X]$. Then, for $h_i =
h_{i,1} - b_i p_{3-i} s$, we have
\[
a = p_i h_i + b_i c \qquad i=1,2.
\]
To conclude the proof it suffices to find $s\in F[X]$ such that
$h_1$ and $h_2$ are separable, $\gcd(h_i, ap_i) = 1$, and $\deg c =
n$. Choose $s\in F[X]$ for which $\deg s = n - (\deg p_1 + \deg
p_2)\geq 1$, $(b p_{i} s)'\neq 0$, and $\gcd(s,h_{i,1}) = 1$ for
$i=1,2$ (e.g., $s(X) = (X-\beta)^{n-1}(X-\gamma)$, where
$\beta,\gamma\in F$ are not roots of $h_{1,1} h_{2,1} b p_1 p_2$.)

By Lemma~\ref{lem_separable} with $h_{i,1}$, $b p_{3-i}s$, $a p_i$
(for $i=1,2$) we get a finite set $S\subseteq F$ such that for each
$\alpha \in F \smallsetminus S$ the polynomial $h_{i,1} - \alpha b
p_{3-i} s$ is separable and relatively prime to $a p_i$. Replace $s$
with $\alpha s$, for some $\alpha\neq 0$ for which
$\alpha_i\alpha\in F\smallsetminus S$, if necessary, to assume that
$\alpha_1,\alpha_2\in F\smallsetminus S$. This $s$ has all the
required properties.
\end{proof}

The next lemma gives a criterion, which we shall use to prove
Proposition~\ref{thmB}, for a transitive group to be primitive, and
further, to be the symmetric group (cf.\ \cite[Lemma
4.4.3]{Serre1992}).

\begin{lemma}\label{lem_primitive}
Let $A\leq S_n$ be a transitive group and let $e$ be a positive
integer in the segment $\frac n2<e<n$ such that $\gcd(e,n)=1$. Then,
if $A$ contains an $e$-cycle, it is primitive. Moreover, if $A$ also
contains a transposition, then $A = S_n$.
\end{lemma}

\begin{proof}
Let $\Delta \neq \{1,\ldots,n\}$ be a block of $A$. We have
$|\Delta| \leq \frac n2$, since $|\Delta|\mid n$. For the first
assertion, it suffices to show that $|\Delta| = 1$, and since
$\gcd(e,n)=1$, it even suffices to prove that $|\Delta|\mid e$.
Without loss of generality assume that $\sigma = (1\ 2\ \cdots\  e)
\in A$ and $1\in \Delta$. Then $\{1,\ldots,e\} \not\subseteq
\Delta$, since $e>\frac n2\geq |\Delta|$. Hence $\Delta \neq \sigma
\Delta$ which implies that $\Delta\cap \sigma \Delta = \emptyset$.
As $\sigma(x) = x$ for any $n\geq x>e$, we have $\Delta \subseteq
\{1,\ldots,e\}$. Consequently, $\Delta$ is a block of
$\left<\sigma\right>$, so $|\Delta|\mid e$.

The second assertion follows since a primitive group containing a
transposition is the symmetric group \cite[Theorem
3.3A]{DixonMortimer}.
\end{proof}

The following number-theoretic lemma will be needed later.

\begin{lemma}\label{lem:numbertheoretic}
For any prime $p$ and positive integers $n,m$ satisfying $n\geq 2m +
\log n (1+o(1))$, there exists an integer $e$ in the segment $\frac
n2< e< n-m$ such that $\gcd(e,np)=1$.
\end{lemma}

\begin{proof}
Let $e$ be
\[
\begin{array}{ll}
\frac{n}{2}+2,
    & \mbox{if $n$ is even but not divisible by $4$,}\\
\frac{n}{2}+1,
    & \mbox{if $n$ is divisible by $4$, or}\\
\frac{n+1}{2},
    & \mbox{if $n$ is odd.}
\end{array}
\]
Then $e$ is the first integer greater than $\frac{n}{2}$ for which
$\gcd(e,n)=1$. If $p\nmid e$, we are done (and we only need $n>2m +
4$). Next assume that $p\mid e$ (and hence $p\nmid n$). Firstly, if
$n$ is even but not divisible by $4$, then the next candidate $e' =
e+2$ works, since $\gcd(e',n)=1$ and $p\nmid e'$. (Otherwise, $p$
divides $e'-e=2$, hence $e$ is even, a contradiction.) Next, if $n$
is divisible by $4$, then the first relatively prime to $n$ integer
greater than $e$ is $e' = \frac{n}{2}+q$, where $q$ is the smallest
prime not dividing $n$. Had $p\mid e'$, we would have $p\mid (e'-e)
= q-1$. In particular, $p<q$, and hence $p\mid n$ by minimality of
$q$, a contradiction. Finally, if $n$ is odd the same argument will
show that $e' = \frac{n+q}{2}$ is relatively prime to $np$, where
now $q$ is the smallest odd prime not dividing $n$ (if $p=2$ we have
to take $q\equiv n \pmod 4$).

It remains to evaluate $q$ which is a standard exercise in number
theory: Let $\omega(n)$ be the number of distinct prime divisors of
$n$. Then $q$ is no more than the $\omega(n)+2$ prime number. Since
the $k$-th prime equals to $k\log k(1+o(1))$ and
$$\omega(n)\leq \frac{\log n}{\log\log n}(1 + o(1))$$
\cite[Theorem 2.10]{MontgomeryVaughan2007}, we have
\[
q\leq \omega(n)\log(\omega(n))(1+o(1)) = \log n (1+o(1)).
\]
Note that for $n = 4 \prod_{2<l<q} l$ (i.e., $4$ times the product
of all the odd prime numbers less than $q$) the inequality is in
fact equality. Thus the estimation $n>2m+\log(n)(1+o(1))$ is the
best possible.
\end{proof}

The following result is very well known, however, for the sake of
completeness, we give a proof.

\begin{proposition} \label{prop:cyclicelements}
Let $F$ be an algebraically closed field of characteristic $l\geq
0$. Let $E/K$ be a separable extension of degree $n$ of algebraic
function fields of one variable over $F$. Assume that a prime
divisor $\frakp$ on $K$ decomposes as
\[
\frakp = \frakP_1^{e_1} \cdots \frakP_r^{e_r}
\]
on $E$. If $l>0$, assume further that $\gcd(e_i,l)=1$, for
$i=1,\ldots,r$. Then the Galois group of the Galois closure of $E/K$
(as a subgroup of $S_n$) contains an element of cyclic type
$(e_1,\ldots,e_r)$. Moreover, the result holds even if $l=e_r=2$ (we
still assume that $\gcd(e_i,2)=1$ for $i=1,\ldots,r-1$).
\end{proposition}

\begin{proof}
The completion $\Khat$ of $K$ at $\frakp$ is a field of Laurent
series over $F$ \cite[Theorem 2]{Serre1979}, say $\Khat = F((Y))$.
Let $x$ be a primitive element of $E/K$, integral at $\frakp$ and
let $f$ be its irreducible polynomial over $K$. Then $f$ factors
over $F((Y))$ into a product of separable irreducible polynomials $f
= f_1 \cdots f_r$ such that $\deg f_i = e_i$ for each $i=1,\ldots,r$
\cite[II\S3]{Serre1979}.

If either $l=0$, or $l>0$ and $\gcd(e_i,l)=1$, then $F((Y))$ has a
unique extension of degree $e_i$, namely $F((Y^{1/e_i}))$
\cite[IV\S6]{Chevalley}. We thus get that the splitting field of $f$
over $F((Y))$ is $F((Y^{1/e}))$, where $e=\lcm(e_1,\ldots,e_r)$,
unless $l=e_r=2$ and then the splitting field of $f$ is the
compositum of $F((Y^{1/e'}))$ with an extension of degree $2$, where
$e' = \lcm(e_1,\ldots,e_{r-1})$. In both cases the Galois group of
$f$ over $F((Y))$ is cyclic of order $e$. Its generator $\sigma$
acts cyclicly on the roots of each of the $f_i$'s. Consequently, the
cyclic type of $\sigma$ is $(e_1,\ldots,e_r)$, as required.
\end{proof}

\begin{lemma}\label{lem:decomposition}
Let $F$ be a field, let $f(X,Y)=a(X) + b(X) Y = \sum_{i=0}^n(a_i
+b_iY)X^i\in F[X,Y]$ be irreducible and separable over $F(Y)$, let
$E=F(Y)[X]/(f(X,Y))$, and let $\alpha\in F$ such that $a_n + b_n
\alpha\neq 0$. Then the decomposition of $\frakp=(Y-\alpha)$ on $E$
corresponds to the factorization of $f(X,\alpha)$ over $F$.
\end{lemma}

\begin{proof}
Let $R = F[Y,(a_n+b_nY)^{-1}]$. Then $S=R[X]/(f(X,Y))\subseteq E$ is
integral over $R$. Moreover, as $Y = -a(X)b^{-1}(X)$ in $E$, we have
that $S$ is a localization of the polynomial ring $F[X]$ at
$-b(X)a^{-1}(X)$ and at $a_n - b_n a(X)b^{-1}(X)$.
Hence $S$ is integrally closed \cite[Proposition
5.12]{AtiyahMacdonald}, and the assertion follows from \cite[I\S4
Proposition 10]{Serre1979}.
\end{proof}

\subsection{Proof of Proposition~\ref{thmB}}
\begin{trivlist}
\item
Let $f(X,Y) = a(X) + b(X)Y\in F[X,Y]$ be an irreducible polynomial.
For large integer $n$ we need to find $c(X)\in F[X]$ such that
$f(X,c(X)Y)=a(X) + b(X) c(X)Y$ is irreducible of degree $n$ and the
Galois group of $f(X,c(X)Y)$ over $\Fgal(Y)$ is $S_n$.

Lemma~\ref{lem:numbertheoretic} with $m=\max\{\deg a(X), 2 + \deg
b(X)\}$ and $p = \Char(F)$ gives (for $n > 2m+\log n(1+o(1))$) a
positive integer $e$ such that
\begin{eqnarray}
&& \label{eq_n_large} n -m > e > \frac{n}{2}\ (\mbox{in particular, } e>m), \\
\label{eq_p_notdivide} &&\gcd(e,np)=1\ (\mbox{or $\gcd(e,n)=1$, if
$p=0$}).
\end{eqnarray}

Let $\alpha_1\neq \alpha_2$ and $\gamma_1\neq \gamma_2$ be elements
of $F$ such that $\alpha_i$ is nonzero and $\gamma_i$ is not a root
of $a(X)b(X)$, $i=1,2$. In Lemma~\ref{lem_final_poly} we constructed
(for $a$, $b$, $p_1 = (X-\gamma_1)^e$, $p_2 = (X-\gamma_2)^2$,
$\alpha_1 $, and $\alpha_2$) a polynomial $c(X)\in F[X]$ of degree
$\deg c = n - \deg b(X)$ which is relatively prime to $a(X)$ such
that
\begin{eqnarray}
\label{eq:ram 1}%
f(X,c(X)\alpha_1 ) &=& a(X) + \alpha_1b(X) c(X) =  (X - \gamma_1)^e h_1(X),\\
\label{eq:ram alpha}%
f(X,c(X)\alpha_2 ) &=& a(X) + \alpha_2 b(X) c(X) = (X-\gamma_2)^2
h_2(X).
\end{eqnarray}
Here $h_1(X),h_2(X)\in F[X]$ are separable polynomials which are
relatively prime to $(X-\gamma_1)a(X)$, $(X-\gamma_2)a(X)$,
respectively. In particular $\gcd(a,c)=1$, and hence $f(X,c(X)Y)$ is
irreducible (Lemma~\ref{lem1}). By \eqref{eq_n_large}, $\deg_X
f(X,c(X) Y)  = \deg b(X) +\deg c(X) = n$. Taking
\eqref{eq_p_notdivide}, \eqref{eq:ram 1}, and \eqref{eq:ram alpha}
in mind, Lemma~\ref{lem:decomposition} and
Proposition~\ref{prop:cyclicelements} with $\frakp = (Y-\alpha_1)$
give us an $e$-cycle in $\gal(f(X,c(X)Y),\Fgal(Y))$ and with
$\frakp=(Y-\alpha_2)$ give a transposition. Thus
$\gal(f(X,c(X)Y),\Fgal(Y))=S_n$ (Lemma~\ref{lem_primitive}). \qed
\end{trivlist}

\section{Dirichlet's theorems}\label{sec:3}
\subsection{Proof of Theorem~\ref{thm A}}
We actually prove a stronger statement:
\begin{theorem}\label{thmalmostA}
Let $M/F$ be a PAC extension. Then Dirichlet's theorem holds for
$F[X]$ and $\calN(M)$.
\end{theorem}

\begin{proof}
The field $F$ is an infinite field, since $M/F$ is a PAC extension
(\cite[Remark 1.2]{JardenRazon1994}), and $f(X,Y) = a(X) + b(X) Y$
is irreducible (Lemma~\ref{lem1}). Proposition~\ref{thmB} gives a
polynomial $c(X)\in F[X]$ for which $f(X,c(X)Y)$ is an irreducible
polynomial of degree $n$ in $X$ and $\gal(f(X,c(X)Y) , \Fgal(Y)) =
S_n$. Now the assertion follows from Corollary~\ref{cor:places}.
\end{proof}

\begin{remark}
The above proves a  stronger statement then stated, namely, for
large $n$ there exists $c(X)\in F[X]$ such that every separable
extension $N/M$ of degree $n$ is generated by a root of
$f(X,c(X)\alpha)$, for infinitely many $\alpha\in F$.
\end{remark}

In Theorem~\ref{thm A} and Proposition~\ref{thmB} we have
considered only linear polynomials. We pose the natural
generalization to general polynomials (cf. \cite[Lemma
10.3.1]{FriedJarden2005}):

\begin{problem}
Let $f(X,Y)\in F[X,Y]$ be an absolutely irreducible polynomial. For
large $n$, is there a polynomial $c(X)\in F[X]$ for which
$f(X,c(X)Y)$ is an $X$-stable polynomial of degree $n$? for which
$\gal(f(X,c(X)Y), \Fgal(Y)) \cong S_n$?
\end{problem}

\subsection{PAC extensions} \label{sec:PAC}
Theorem~\ref{thm A} makes it interesting to calculate $\calN(F)$ for
a PAC field $F$. Its generalization, Theorem~\ref{thmalmostA},
raises the following questions. When does a given field have a PAC
extension $M$? and what positive integers can occur as degrees of
separable extensions of such $M$'s?

\begin{examples}[PAC fields which have separable extensions of arbitrary
degrees] Let $M$ be a PAC field. If $M$ is Hilbertian or more
generally RG-Hilbertian (i.e., $M$ has the irreducible
specialization property for regular Galois extensions), then every
finite group occurs as a Galois group over $M$
\cite{FriedVolklein1992}. In particular, $\calN(M)=\mathbb Z^+$
(where $\mathbb Z^+$ denotes the set of all positive integers).

In general, $n\in \calN(M)$ if and only if some Galois group over
$M$ has a subgroup of index $n$. For example, if a cyclic group of
order $n$ (or alternatively the symmetric group of degree $n$)
occurs as a Galois group over $M$, then $n\in\calN(M)$. In
particular, if $\gal(M)$ is a finitely generated free profinite
group (and hence $M$ is ``far'' from being Hilbertian), then
$\calN(M)=\mathbb Z^+$.
\end{examples}

The succeeding result asserts that for a PAC field $M$, the set
$\calN(M)$ is finite only if $M= M_s$, where $M_s$ is a fixed
separable closure of $M$.

\begin{lemma}\label{lem:infinite index}
Let $M$ be a PAC field and assume $M\neq M_s$. Then $\calN(M)$ is
infinite.
\end{lemma}

\begin{proof}
Artin-Schreier Theorem implies that if $[M_s:M]$ is finite, then
$[M_s:M]=2$ and $M$ is real closed. However a PAC field cannot be
real, since the curve defined by $X^2 + Y^2 + 1$ has an $M$-rational
point.
\end{proof}

In contrast to PAC fields which have been studied since the late
1960s, PAC extensions first appeared in 1994 and have not yet been
well understood. In what follows we describe all the known PAC
extensions. Also we give some new explicit interesting examples of
PAC extensions and deduce for them corollaries to
Theorem~\ref{thmalmostA}.

A profinite group is compact, hence it is equipped with a normalized
Haar measure (see \cite[Chapter 18]{FriedJarden2005}). In
particular, if $K$ is a field and $e$ is a positive integer, then
$\gal(K)^e$ is equipped with a normalized Haar measure. Let
$\bfsigma = (\sigma_1\nek \sigma_e) \in \gal(K)^e$ be an $e$-tuple
of Galois automorphisms. Then $\left< \bfsigma\right>$ denotes the
subgroup generated by $\sigma_1\nek \sigma_e$ and $K_s(\bfsigma)$
denotes the fixed field of $\left<\bfsigma\right>$ in a fixed
separable closure $K_s$ of $K$. The phrase ``for almost all'' means
``for all but a set of measure zero.''

Clearly (and uninteresting for us), a separably closed field is PAC
extension of any infinite subfield of it. In \cite{JardenRazon1994},
Razon and Jarden prove the following:
\begin{enumerate}
\item \label{foralmostall}
Let $K$ be a countable Hilbertian field and let $e\geq 1$. Then
$K_s(\bfsigma)/K$ is PAC (as an extension) for almost all
$\bfsigma\in \gal(K)^e$.
\item \label{PACextensionextension}
Let $L/K$ be a PAC extension and let $F/K$ be an algebraic
extension. Then $FL/F$ is a PAC extension.
\end{enumerate}
So far these are all the known PAC extensions. We shall use these
properties to construct some explicit PAC extensions $M/F$ (for
non-Hilbertian $F$). Before doing that, we emphasize that there are
restrictions on a field extension to be PAC, e.g., it is known that
no Galois extension of a finitely generated field (except for the
separable closure) is PAC \cite{Bary-SorokerJarden}. In particular,
since there is a Galois extension $M$ of $\bbQ$ such that $M$ is a
PAC field \cite[Theorem 18.10.2]{FriedJarden2005}, we get that for
an extension $M/F$ to be PAC (as an extension) it does not suffice
that $M$ is PAC (as a field).

Given a countably Hilbertian field $K$, \eqref{foralmostall} and
\eqref{PACextensionextension} yield abundance of separable
extensions $F/K$ such that Dirichlet's theorem holds for $F[X]$ (and
some infinite set $\calN$ which we omit from now on):

\begin{corollary}
\label{cor:DTpositivemeasure} Let $K$ be a countable Hilbertian
field and let $F/K$ be a separable algebraic extension. Assume that
the set $\{\bfsigma \in \gal(K)^e\mid FK_s(\bfsigma) \neq K_s\}$ has
a positive measure for some positive integer $e$. Then, Dirichlet's
theorem holds for $F[X]$.
\end{corollary}

\begin{proof}
By the assumption, \eqref{foralmostall}, and
\eqref{PACextensionextension} we get that there is $\bfsigma$
(actually a positive measure set of $\bfsigma$'s) for which the
extension $FK_s(\bfsigma)/F$ is PAC and $FK_s(\bfsigma)\neq K_s$.
Therefore Lemma~\ref{lem:infinite index} and
Theorem~\ref{thmalmostA} imply the assertion.
\end{proof}

For almost all $\bfsigma\in \gal(K)^e$ the group
$\left<\bfsigma\right>$ is isomorphic to the free profinite group on
$e$ generators \cite[Theorem 18.5.6]{FriedJarden2005}. Also, by
Galois correspondence, if $F/K$ is Galois and $FK_s(\bfsigma) =
K_s$, then $\left<\bfsigma\right>$ is isomorphic to a subgroup of
$\gal(F/K)$. Therefore we have

\begin{corollary}
\label{prop4}%
Let $K$ be a countable Hilbertian field, let $e\geq 1$, let $F/K$ be
a Galois extension, and assume that $\gal(F/K)$ does not have a free
subgroup of rank $e$. Then Dirichlet's theorem holds for $F[X]$.
\end{corollary}

In particular, since a pro-solvable group cannot have a noncyclic
free profinite group as a subgroup, we get

\begin{corollary}
Let $F$ be a pro-solvable extension of a countable Hilbertian
field $K$. Then Dirichlet's theorem holds for $F[X]$.
\end{corollary}

In light of the above corollaries we suggest two open problems:

\begin{problem}
Let $K$ be a Hilbertian field. Classify all separable algebraic
extensions $F$ of $K$ such that the measure of $\{ \bfsigma\in
\gal(K)^e \mid F K_s(\bfsigma) \neq K_s\}$ is positive for some
positive integer $e$.
\end{problem}

\begin{problem}
Let $K$ be a Hilbertian field. Classify extensions $F$ of $K$ which
have a PAC extension which is not separably closed.
\end{problem}

\begin{remark}
There are fields which do not have any separable extension which
is PAC, other than the separable closure, e.g. $\bbC$, $\bbR$,
$\bbQ_p$. In general, if a field is henselian, then it has no
separable extension which is a PAC field other than the separable
closure \cite[Corollary 11.5.5]{FriedJarden2005}. Note that over
$\bbC$ or over $\bbR$ Dirichlet's theorem obviously does not hold,
since every polynomial of degree greater than two is reducible.
\end{remark}


\begin{thebibliography}{MMM99}

\bibitem[AM69]{AtiyahMacdonald}
Michael F. Atiyah and Ian. G. Macdonald, \emph{Introduction to
commutative algebra}, Addison-Wesley Publishing Co., Reading,
Mass.-London-Don Mills, Ont., 1969.

\bibitem[BSJ]{Bary-SorokerJarden}
Lior Bary-Soroker and Moshe Jarden, \emph{{PAC} fields over finitely
generated fields}, to appear in Math. Z.

\bibitem[Che51]{Chevalley}
Claude Chevalley, \emph{Introduction to the theory of algebraic
functions of one variable}, Mathematical Surveys, no. 6, American
Mathematical Society, New York, 1951.

\bibitem[DM96]{DixonMortimer}
John~D. Dixon and Brian Mortimer, \emph{Permutation groups},
Graduate Texts in
  Mathematics, vol. 163, Springer-Verlag, New York, 1996.

\bibitem[FHJ84]{FriedHaranJarden1984}
Michael~D. Fried, Dan Haran, and Moshe Jarden \emph{Galois
stratification over Frobenius fields}, Adv. in Math. \textbf{51}
(1984), no. 1, 1--35.

\bibitem[FJ05]{FriedJarden2005}
Michael~D. Fried and Moshe Jarden, \emph{Field arithmetic}, second
ed., revised and enlarged by Moshe Jarden,
  Ergebnisse der Mathematik (3) \textbf{11}, Springer-Verlag, Heidelberg, 2005.

\bibitem[FV92]{FriedVolklein1992}
Michael~D. Fried and Helmut V\"olklein,  \emph{The embedding problem
over a Hilbertian PAC-field}, Ann. of Math. (2) \textbf{135} (1992),
no. 3, 469--481.


\bibitem[JR94]{JardenRazon1994}
Moshe Jarden and Aharon Razon, \emph{Pseudo algebraically closed
fields over
  rings}, Israel J. Math. \textbf{86} (1994), no.~1-3, 25--59.

\bibitem[Kor19]{Kornblum1919}
Heinrich Kornblum, \emph{\"{U}ber die {P}rimfunktionen in einer
arithmetischen {P}rogression}, Mathematische Zeitschrift \textbf{5}
(1919), 100--111.

\bibitem[MM99]{MalleMatzat1999}
Gunter Malle and B.~Heinrich Matzat, \emph{Inverse {G}alois theory},
Springer Monographs in Mathematics, Springer-Verlag, Berlin, 1999.

\bibitem[MV07]{MontgomeryVaughan2007}
Hugh L.~Montgomery and Robert C.~Vaughan, \emph{Multiplicative
Number Theory I}, Cambridge Studies in Advance Mathematics
\textbf{97}, Cambridge University Press, Cambridge, 2007.


\bibitem[Ros02]{Rosen2002}
Michael Rosen, \emph{Number theory in function fields}, Graduate
Texts in
  Mathematics, \textbf{210}, Springer-Verlag, New York, 2002.

\bibitem[Ser79]{Serre1979}
Jean-Pierre Serre, \emph{Local fields}, Graduate Texts in
Mathematics, \textbf{67}, Springer-Verlag, New York, 1979.



\bibitem[Ser92]{Serre1992}
Jean-Pierre Serre, \emph{Topics in {G}alois theory}, Research Notes
in Mathematics, \textbf{1}, Jones and Bartlett Publishers, Boston,
1992.
\end{thebibliography}
\end{document}